\documentclass[letterpaper, 10 pt, conference]{cssconf}  
\IEEEoverridecommandlockouts                              
\usepackage{graphicx}      

\usepackage{amsmath}
\usepackage{amssymb}

\usepackage{url}

\usepackage{tikz}

\newcommand{\BBM}{\begin{bmatrix}}
\newcommand{\EBM}{\end{bmatrix}}
\newcommand{\BEAS}{\begin{eqnarray*}}
\newcommand{\EEAS}{\end{eqnarray*}}
\newcommand{\BEA}{\begin{eqnarray}}
\newcommand{\EEA}{\end{eqnarray}}
\newcommand{\BEQ}{\begin{equation}}
\newcommand{\EEQ}{\end{equation}}
\newcommand{\BIT}{\begin{itemize}}
\newcommand{\EIT}{\end{itemize}}

\newcommand{\reals}{{\mbox{\bf R}}}
\newcommand{\naturals}{{\mbox{\bf N}}}

\newcommand{\diag}{\mathop{\bf diag}}

\title{\LARGE \bf Efficient Robust Model Predictive Control using Chordality}

\author{Anders Hansson and Sina Khoshfetrat Pakazad 
\thanks{This research has been supported by ELLIIT and by 
the Wallenberg Artificial Intelligence, Autonomous Systems and Software 
Program (WASP) funded by Knut and Alice Wallenberg Foundation, 
which is gratefully acknowledged. The authors are also grateful for
discussion with Dimitris Kouzoupis. }
\thanks{Division of Automatic Control, Link\"oping University, Sweden, 
 {\tt\small anders.g.hansson@liu.se}}%
\thanks{C3 IoT, Redwood city, California, USA, 
 {\tt\small sina.pakazad@c3iot.com}}%
}

\begin{document}

\maketitle
\thispagestyle{empty}
\pagestyle{empty}

\begin{abstract}                
In this paper we show that chordal structure can be used
to devise efficient optimization methods for robust model predictive
control problems. The chordal structure is used both for 
computing search directions efficiently as well as for distributing all the
other computations in an interior-point method for solving the problem. 
The framework enables efficient parallel computations.
\end{abstract}

\section{Introduction}
Model Predictive Control (MPC) is an important class of controllers that 
are being employed more and more in industry, 
\cite{qin+bad03}. It has its root
going back to \cite{cut+ram79}. The success is mainly because it can handle 
constraints on control signals and/or states in a systematic way. In the
early years its applicability was limited to slow processes, since 
an optimization problem has to be solved at each sampling instant. 
Tremendous amount of research has been spent on overcoming this limitation. 
One avenue has been what is called explicit MPC, \cite{Alessio2009}, 
where the optimization
problem is solved parametrically off-line. Another avenue has been to 
exploit the inherent structure of the optimization problems stemming from MPC, 
\cite{gla+jon84,wri93,ste94,arn+put94,wri96,rao+wri+raw97,gop+bie98,han98,van+boy+nou01,jor04,ake+han04,Diehl2009,axe+van+han10,wan+boy10,jer+ker+con12,dom+12,fri15,kli17,nie17}. Typically
this has been to use Riccati recursions to efficiently compute search directions for 
Interior Point (IP) methods or actives set methods to solve the optimization problem.
In \cite{han+pak17} it was argued that the important structures that have been 
exploited can all be summarized as {\it chordal structure}. Because of this
the same structure exploiting software can be used to speed up all computations
for MPC. This is irrespective of what MPC formulation is considered and 
irrespective of what type of optimization algorithm is used. In this paper we will 
in detail discuss robust MPC, which was not discussed in the above mentioned reference. 
We assume that the reader is familiar with the receding horizon strategy
of MPC and we will only discuss the associated constrained finite-time
optimal control problem. We will from now on refer to 
the associated problem as the MPC problem.

The remaining part of the paper is organized as follows.
We will in Section~2 discuss how chordal sparsity arises and how it can be 
utilized in general convex optimization problems to obtain 
computations distributed over a so called clique tree. The presentation is based
on \cite{khoshfetrat2016distributed}. In Section~3 we then discuss how this
is can be used within IP methods for general Robust Quadratic Programs (RQPs). 
In Section~4 we state the robust  MPC problem. It is formulated using a 
scenario tree, and we will see that it is a special case of RQP. 
In Section~5 we will give some conclusions, 
discuss generalizations of our results 
and directions for future research. 
\subsection*{Notation}
We denote with $\reals$ the set of real numbers, with $\reals^n$ the set of 
$n$-dimensional real-valued vectors and with $\reals^{m\times n}$ the set of 
real-valued matrices with $m$ rows and $n$ columns. 
We denote by $\naturals$ the set of natural
numbers and by $\naturals_n$ the subset  $\lbrace 1,2,\ldots,n\rbrace$ of $\naturals$.
For a vector $x\in\reals^n$ the matrix $X=\diag(x)$ is a diagonal matrix with the
components of $x$ on the diagonal. For two matrices $A$ and $B$ the matrix $A\oplus B$ is a block-diagonal matrix with 
$A$ as the 1,1-block and $B$ as the 2,2-block. For a symmetric matrix $A$ 
the notation $A(\succeq)\succ 0$ is equivalent to $A$ being positive (semi)-definite. 
\section{Chordal Sparsity and Convex Optimization}
Consider the following convex optimization problem
\begin{align}\label{eqPaperIII:CP}
\min_{x} \quad F_1(x) + \dots + F_N(x),
\end{align}
where $F_i \ : \ \reals^{n}\rightarrow \reals$ for all $i = 1, \dots, N$. We assume that each function $F_i$ is only dependent on a small subset of elements of $x$. Let us denote the ordered set of these indexes by $J_i \subseteq \naturals_n$. We can then rewrite the problem in~\eqref{eqPaperIII:CP}, as
\begin{align}\label{eqPaperIII:CPS}
\min_{x} & \quad   \bar F_1(E_{J_1}x) + \dots + \bar F_N(E_{J_N}x),
\end{align}
where $E_{J_i}$ is a $0$--$1$ matrix that is obtained from an identity matrix of order $n$ by deleting the rows indexed by $\naturals_n \setminus J_i$.  The functions $\bar F_i \ : \ \reals^{|J_i|} \rightarrow \reals$ are lower dimensional descriptions of $F_i$s such that $F_i(x) = \bar F_i(E_{J_i}x)$ for all $x \in \reals^n$ and $i \in\naturals_N$. For details on how this structure can be exploited using message passing the reader is referred to \cite{khoshfetrat2016distributed}. 

A brief summary is that we may define a so-called sparsity graph for the above optimization problem with $n$ nodes and edges between two nodes $j$ and $k$ if $x_j$ and $x_k$ appear in the same term $\bar F_i$. We assume that this graph is chordal, i.e. every cycle of length four our more has a chord.\footnote{In case the graph is not chordal we make a chordal embedding, i.e. we add edges to the graph until it becomes chordal. This corresponds to saying that some of the $\bar F_i$ depend on variables that they do not depend on.} The maximal complete subgraphs of a graph are called its cliques. If the original graph is chordal then there exists a tree of the cliques called the clique tree which is such that it enjoys the clique intersection property. This property is that all elements in the intersection of two cliques $C_i$ and $C_j$ should be elements of the cliques on the path between the cliques $C_i$ and $C_j$. It is then possible to use the clique tree as a computational tree where we non-uniquely assign terms of the objective function to each clique in such a way that all the variables of the term in the function are elements of the clique. After this we may solve the optimization problem distributively over the clique tree by starting with leafs and for each leaf solve a parametric optimization problem, where we optimize with the respect to the variables of the leaf problem which are not variables of the parent of the leaf in the clique tree. The optimization should be done parametrically with respect to all the variables that are shared with the parent. After this the optimal objective function value of the leaf can be expressed as a function of the variables that are shared with the parent. This function is sent to the parent and added to its objective function term. The leaf has been pruned away, and then the optimization can continue with the parent assuming all its children has also carried out their local optimizations. Eventually we reach the root of the tree, where the remaining variables are optimized. Then we can finally go down the tree and recover all optimal variables. This is based on the fact that we have stored the parametric optimal solutions in the nodes of the clique tree.  
\section{Interior-Point Methods}
The robust MPC problem is a special case of a so-called Robust 
Quadratic Program (RQP). We will now discuss
how such a problem can be solved using IP methods, \cite{wri97}. 
Consider the RQP
\begin{subequations}\label{eq:RQP}
\begin{align}
\min_{\tau,t,z} \;&\tau\label{t:objective}\\
{\rm s.t.}\;&\frac{1}{2}\left(z_0^j\right)^T\mathcal Q_0^j z_0^j+t_1^j\leq \tau,\quad j\in \naturals_M
\label{RQP:objective0}\\
&\frac{1}{2}\left(z_k^j\right)^T\mathcal Q_k^j z_k^j+t_{k+1}^j\leq t_k^j,\quad 
j\in \naturals_M,\;k\in \naturals_{N-1}
\label{RQP:objectivek}\\
&\frac{1}{2}\left(z_N^j\right)^T\mathcal Q_N^j z_N^j\leq t_N^j,\quad 
j\in \naturals_M,\label{RQP:objectiveN}\\
& \mathcal A z = b\label{RQP:eq}\\
&\mathcal C z\leq d\label{RQP:ineq}
\end{align}
\end{subequations}
where $\mathcal Q_k^j\succeq 0$, i.e. positive semidefinite, where 
$\mathcal A$ has 
full row rank, and where the matrices and vectors are of compatible dimensions.
Here $z=(z^1,\ldots,z^M)$ with $z^j=(z_0^j,\ldots,z_N^j)$, 
and the inequality in (\ref{RQP:ineq})
is component-wise inequality. We will detail the dimensions of $z_k^j$, $b$ and
$d$ later on. Introduce 
$t=(t^1,\ldots,t^M)$ with $t^j=(t_1^j,\ldots,t_N^j)$. We let 
\begin{align*}
\mathcal Q^j&=\oplus_{k=0}^N\mathcal Q_k^j\\
\mathcal Q&=\oplus_{j=1}^M\mathcal Q^j\\
\mathcal Q_\mu^j&=\oplus_{k=0}^N\mu_k^j\mathcal Q_k^j\\
\mathcal Q_\mu&=\oplus_{j=1}^M\mathcal Q_\mu^j\\
\mathcal Q_z^j&=\oplus_{k=0}^N\mathcal Q_k^jz_k^j\\
\mathcal Q_z&=\oplus_{j=1}^M\mathcal Q_z^j
\end{align*}
where $\mu_k^j\geq 0$ are Lagrange
multipliers for the inequality constraints in 
(\ref{RQP:objective0}--\ref{RQP:objectiveN}). We also define
$$\mathcal B^j=\begin{bmatrix}1  &  &       &       &\\
                              -1 & 1&       &       &\\
                                 &-1& 1     &       &\\
                                 &  &\ddots &\ddots &\\
                                 &  &       &     -1& 1\\
                                 &  &       &       &-1\end{bmatrix}\in
\reals^{(N+1)\times N}$$
and $\mathcal B=\oplus_{j=1}^M\mathcal B^j$. Finally we let 
$\beta=(e_1,\ldots,e_1)$, were $e_1\in\reals^{N+1}$ is the first unit vector, 
and where $\beta$ contains $M$ of these vectors. The Lagrangian for the 
optimization problem may now be written as
\begin{align*}
L(\tau,t,z,\mu,\nu\lambda)&=(1-\mu^T\beta)\tau-\mu^T\mathcal B t+
\frac{1}{2}z^T\mathcal Q_\mu z\\
&+\nu^T(\mathcal Cz-d)+\lambda^T(\mathcal Az-b)
\end{align*}
where $\lambda$ and $\nu$ are Lagrange multipliers for the constraints in
\eqref{RQP:eq} and \eqref{RQP:ineq}, respectively, and where 
$\mu=(\mu^1,\ldots,\mu^M)$ with $\mu^j=(\mu_0^j,\ldots,\mu_N^j)$ are the multipliers associated to the remaining constraints.
The Karush-Kuhn-Tucker (KKT) optimality conditions for this problem are
\begin{align*}
1-\beta^T\mu&=0\\
\mathcal B^T\mu&=0\\
\mathcal Q_\mu z+\mathcal A^T\lambda+\mathcal C^T\nu&=0\\
\frac{1}{2}\mathcal Q_z^T z-\beta\tau+\mathcal B t+s&= 0\\
\mathcal A z &= b\\
\mathcal C z+w&= d\\
\mu_k^js_k^j&=0\\
\nu_k^jw_k^j&=0\\
\label{eqn:KKT_RQP}
\end{align*}
and $(\mu,\nu,s,w)\geq 0$, where the vectors $s$ and $w$ are slack variables
for the inequality constraints.

In IP methods one linearizes
the above equations to obtain equations for search directions
\begin{align*}
-\beta^T\Delta\mu&=r_\tau\\
\mathcal B^T\Delta\mu&=r_t\\
\mathcal Q_\mu \Delta z+\mathcal Q_z\Delta\mu+\mathcal A^T\Delta\lambda+
\mathcal C^T\Delta\nu&=r_\mu\\
\mathcal Q_z^T\Delta z+\mathcal B\Delta t-\beta\Delta\tau+\Delta s&=r_z\\
\mathcal A\Delta z&=r_\lambda\\
\mathcal C\Delta z+\Delta w&=r_\nu\\
M\Delta s+S\Delta\mu&=r_s\\
V\Delta w+W\Delta\nu&=r_w
\end{align*}
where $M=\diag(\mu)$, $S=\diag(s)$, $V=\diag(\nu)$, $W=\diag(w)$, 
and where $r=(r_t,r_\mu,r_z,r_\lambda,r_\nu,r_s,r_w)$ is some
residual vector that depends on what IP method is used. The
quantities $r$, $M$, $S$, $V$, $W$, $\mathcal Q_\mu$ and $\mathcal Q_z$
depend on the value of the current 
iterate in the IP method.

From the last three
equations above we have $\Delta w=r_\nu-\mathcal C\Delta z$, 
$\Delta s=M^{-1}(r_s-S\Delta \mu)$ and $\Delta\nu
=W^{-1}(r_w-V\Delta w)$. After substitution of these expressions into the third
and fourth equation we obtain
\begin{align*}
(\mathcal Q_\mu+\mathcal C^TW^{-1}V\mathcal C)\Delta z+
\mathcal Q_z\Delta\mu+\mathcal A^T\Delta\lambda&=\bar r_\mu\\
\mathcal Q_z^T\Delta z+\mathcal B\Delta t-\beta\Delta\tau-M^{-1}S\Delta\mu
&=\bar r_z
\end{align*}
where $\bar r_\mu=r_\mu-\mathcal C^TW^{-1}(r_w-Vr_\nu)$ and
$\bar r_z=r_z-M^{-1}r_s$.
Solve from the last equation with respect to $\Delta\mu$ to obtain
$$\Delta\mu=S^{-1}M(\mathcal Q_z^T\Delta z+\mathcal B\Delta t-\beta\Delta\tau
-\bar r_z)$$
We now substitute this expression into all equations containing $\Delta\mu$
and obtain the following linear system of equations for the remaining variables
$$\begin{bmatrix}P&R^T&\\
                 R&\mathcal Q_s&\mathcal A^T\\
                  &\mathcal A\end{bmatrix}
\begin{bmatrix}
\Delta\eta\\\Delta z\\\Delta\lambda
\end{bmatrix}=
\begin{bmatrix}
r_\eta\\\tilde r_\mu\\r_\lambda
\end{bmatrix}
$$
where $\Delta\eta=(\Delta\tau,\Delta t)$, 
\begin{align*}
P&=\begin{bmatrix}-\beta&\mathcal B\end{bmatrix}^TS^{-1}M
\begin{bmatrix}-\beta&\mathcal B\end{bmatrix}\\
R&=\mathcal Q_zS^{-1}M
\begin{bmatrix}-\beta&\mathcal B\end{bmatrix}\\
\mathcal Q_s&=\mathcal Q_\mu+\mathcal C^TW^{-1}V\mathcal C +\mathcal Q_z
S^{-1}M\mathcal Q_z^T
\end{align*}
and where $r_\eta=(r_\tau-\beta^TS^{-1}M\bar r_z,r_t-\mathcal B^TS^{-1}M\bar r_z)$.
The matrix $P$ is invertible by construction. 
Hence we may solve for $\Delta \eta$ to obtain 
$\Delta\eta=P^{-1}(r\eta-R^T\Delta z)$ and 
substitute into the other equations, which gives
\begin{equation}\label{eqn:KKT-reduced}
\begin{bmatrix}\mathcal Q_s-RP^{-1}R^T&\mathcal A^T\\
                  \mathcal A&\end{bmatrix}
\begin{bmatrix}
\Delta z\\\Delta\lambda
\end{bmatrix}=
\begin{bmatrix}
\tilde r_\mu-P^{-1}r_\eta\\r_\lambda
\end{bmatrix}
\end{equation}
We notice that the search directions are 
obtained by solving an indefinite symmetric
linear system of equations. This matrix is referred to as the KKT matrix. 
Notice that the above linear system of equations 
for the search directions can be interpreted as the optimality conditions 
for a QP in $\Delta z$ 
with only equality constraints. 
In case this QP is loosely coupled with chordal 
structure message passing over a clique tree can be used to compute
the search directions in a distributed way. 
This is explained in more detail in \cite{khoshfetrat2016distributed,han+pak17}.
We remark that in case $\mathcal A$ and $\mathcal C$ are block diagonal , i.e. 
there is no coupling for different $j$ in the constraints, then
the only coupling with respect to $j$ is related to $RP^{-1}R^T$. The
matrix $P$ has a block arrow structure, and the coupling structure is weak 
because of the structure of $\beta$. It is easy to see that only the variables
$\Delta z_0^j$ will be coupled. We will see that we have even more structure
that can be exploited for robust MPC. 
\section{Robust MPC}
There are many ways to define robust (linear)
MPC problems. However, they all fall
into the category 
\begin{align*}
\min_{u(p)}&\max_{p\in\mathcal P}
\frac{1}{2}\sum_{k=0}^{N-1}\begin{bmatrix}x_k(p)\\u_k(p)\end{bmatrix}^TQ
\begin{bmatrix}x_k(p)\\u_k(p)\end{bmatrix}+\frac{1}{2}x_N(p)^TSx_N(p)\\
{\rm s.t.}\; &x_{k+1}(p)=A(p)x_k(p)+B(p)u_k(p)+v_k(p),\; x_0=\bar x\\
&Cx_k(p)+Du_k(p)\leq e_k
\end{align*}
where $\mathcal P$ is some set. Here $A(p)\in\reals^{n\times n}$ and
$B(p)\in \reals^{n\times m}$. Here $e_k$ is not a basis vector. 
We also assume that there are $q$ inequality 
constraints for each $k$ and that the dimensions of the other matrices and
vectors are compatible with this. One usually makes the assumption that
$p$ depends on $k$ and that $u_k$ only depends on values of $p$ prior to
$k$, the so-called non-anticapativity constraint. 
Since point-wise maximum over convex functions preserves convexity, 
it follows that
the above problem also is convex. It should, however, be stressed that
it is in general not tractable unless further assumptions are made on 
$\mathcal P$, such as e.g. finiteness.
It is possible to also let $C$, $D$, $Q$, $S$, and
$e_k$ depend on $p$ without destroying convexity. 

We will consider a special important case
that is obtained by letting the dynamics evolve as
\begin{multline*}
x_{k+1}(\bar p_k)=A(p_k)x_k(\bar p_{k-1})+\\ B(p_k)u_k(\bar p_{k-1})+v_k(p_k),
\; x_0=\bar x
\end{multline*}
where $\bar p_k=(p_0,p_1\ldots,p_k)$, with $p_k\in\mathcal P_k$, where
$\mathcal P_k$ are finite sets with cardinality $M_k$.
We realize that 
the number of equality constraints
grows exponentially with $k$, in case
the cardinality is independent of $k$. In order to get tractable problems
one often let $M_k=1$ for $k>N_r$, for some integer $N_r$.
Then the problem can be written
\begin{align*}
\min_{u}&\max_{\bar p_{N-1}\in\bar{\mathcal P}_{N-1}}
\frac{1}{2}\sum_{k=0}^{N-1}
\begin{bmatrix}x_k(\bar p_{k-1})\\u_k(\bar p_{k-1})
\end{bmatrix}^TQ
\begin{bmatrix}x_k(\bar p_{k-1})\\u_k(\bar p_{k-1})\end{bmatrix}\\
&+\frac{1}{2}x_N(\bar p_{N-1})^TSx_N(\bar p_{N-1})\\
{\rm s.t.}\;& x_{k+1}(\bar p_k)=
A(p_k)x_k(\bar p_{k-1})+B(p_k)u_k(\bar p_{k-1})+v_k(p_k)\\
&Cx_k(\bar p_{k-1})+Du_k(\bar p_{k-1})\leq e_k
\end{align*}
where $x_0=\bar x$, $u=(u_0,u_1(\bar p_0),\ldots,u_{N-1}(\bar p_{N-1}))$, and
where $\bar{\mathcal P}_k={\mathcal P}_0\times {\mathcal P}_1\times \cdots
\times {\mathcal P}_k$. 

We will now reformulate the problem into an equivalent problem with more
variables and constraints. We let all states and control signals depend
on $p=\bar p_{N_r}\in \mathcal P = \bar{\mathcal P}_{N_r}$ with 
cardinality $M=M_0\times M_1\times\cdots\times M_{N_r}$, i.e. we introduce 
$M$ independent scenarios which we constrain using so-called non-anticipativity
constraints:
\begin{align*}
\min_{u}&\max_{p\in\mathcal P}
\frac{1}{2}\sum_{k=0}^{N-1}
\begin{bmatrix}\bar x_k(p)\\\bar u_k(p)
\end{bmatrix}^TQ
\begin{bmatrix}\bar x_k(p)\\\bar u_k(p)\end{bmatrix}
+\frac{1}{2}x_N(p)^TSx_N(p)\\
{\rm s.t.}\;& \bar x_{k+1}(p)=
A(p_k)\bar x_k(p)+B(p_k)\bar u_k(p)+v_k(p_k)\\
&Cx_k(p)+Du_k(p)\leq e_k
\end{align*}
where $x_0(p)=\bar x$, 
\begin{align*}
&\bar u_k(p_0,\ldots,p_{k},p_{k+1}^1,\ldots,p_{N_r}^1)=\\
&\bar u_k(p_0,\ldots,p_{k},p_{k+1}^2,\ldots,p_{N_r}^2)
\end{align*} 
for all 
$p_{k+1}^1,\ldots,p_{N_r}^1;\;p_{k+1}^2,\ldots,p_{N_r}^2$, and where
$$u=(\bar u_0(p),\ldots,\bar u_{N-1}(p))$$ 
We further define an enumeration of all scenarios using an index $j\in
\lbrace 1,2,\ldots, M\rbrace$ which make it possible to define
the equivalent problem
\begin{align*}
\min_{u}\max_{1\leq j\leq M}&
\frac{1}{2}\sum_{k=0}^{N-1}
\begin{bmatrix}x_k^j\\u_k^j\end{bmatrix}^TQ
\begin{bmatrix}x_k^j\\u_k^j\end{bmatrix}+
\frac{1}{2}(x_N^j)^TSx_N^j\\
{\rm s.t.}\;& x_{k+1}^j=A_k^jx_k^j+B_k^ju_k^j+v_k^j,\quad x_0^j=\bar x\\
&Cx_k^j+Du_k^j\leq e_k\\
&\bar C u = 0
\end{align*}
where $u=(u^1,u^2,\ldots,u^M)$ with $u^j=(u_0^j,u_1^j,\ldots,u_{N-1}^j)$, and  
$$\bar C=\begin{bmatrix}
C_{1,2}&-C_{1,2}&      &        &        \\
      &C_{2,3}&-C_{2,3}&        &        \\
      &      &\ddots &\ddots &         \\
      &      &       &C_{M-1,M}&-C_{M-1,M}
\end{bmatrix}$$ 
with 
$$C_{j,j+1}=\begin{bmatrix}I&0\end{bmatrix}$$
where $I$ is an identity matrix of dimension $m$ times the number of
time instances that scenarios $j$ and $j+1$ have a control signal in 
common. Notice that several
of the matrices $A_k^j$, $B_k^j$ and $v_k^j$ are also constrained, however, we do not have to write those out as they are not optimization variables. Exploiting structure stemming from scenario trees have been investigated in a stochastic setting, e.g. \cite{gon+gro07,mar+15,lei+pot+boc15,fri+17,kou+18}. Here we show how this structure can be exploited due to chordality of the inherent coupling in the problem.

The above problem is equivalent with the problem in \eqref{eq:RQP}. To see this we let $Q_k^j=Q$ and $z_k^j=(x_k^j,u_k^j)$
for $k=0,\ldots,N-1$, and $Q_N^j=S$ and $z_N^j=x_N^j$. We also let
\begin{align*}
b^j&=(\bar x,v_0^j,v_1^j,\ldots,v_{N-1}^j)\\
e^j&=(e_0,e_1,\ldots,e_{N-1})
\end{align*}
and 
\tiny
\begin{align*}
\mathcal A^j&=
\begin{bmatrix}
 I &   &   &   &  &      &   &   &  \\
 -A_0^j&-B_0^j &I  &   &  &      &   &   &  \\
   &   &-A_1^j &-B_1^j &I &      &   &   &  \\
   &   &   &   &  &\ddots&   &   &  \\
   &   &   &   &  &      &-A_{N-1}^j &-B_{N-1}^j &I \\
\end{bmatrix}\\
\mathcal D^j&=\begin{bmatrix}C&D\end{bmatrix}\oplus
\begin{bmatrix}C&D\end{bmatrix}\oplus\ldots\oplus
\begin{bmatrix}C&D\end{bmatrix}
\end{align*}
\normalsize
Finally we let $\mathcal A = \oplus_{j=1}^M\mathcal A^j\oplus\tilde C$, 
$\mathcal D = \oplus_{j=1}^M\mathcal D^j$,  
$b=(b^1,\ldots,b^M)$ and $e=(e^1,\ldots,e^M)$. Here $\tilde C$ is a 
matrix obtained from $\bar C$ by combining its columns with zero columns
such that the non-anticipativity constraint holds. 

We see that the data matrices are banded. Hence, sparse linear system 
solvers could be used  when solving 
\eqref{eqn:KKT-reduced} for 
search directions in an IP method, but we 
will see that the structure within the bands
can be further utilized. 
The only coupling between the $N$ different time instances for a fixed $j$
is via the dynamic equation for $j$. 
The associated QP can be written
\begin{align*}
\min_{\Delta u,\Delta x}&
\begin{bmatrix}\Delta x_0\\\Delta u_0\end{bmatrix}^T\tilde Q_0
\begin{bmatrix}\Delta x_0\\\Delta u_0\end{bmatrix}+\tilde r_0^T
\begin{bmatrix}\Delta x_0\\\Delta u_0\end{bmatrix}\\
&+\frac{1}{2}\sum_{j=1}^M
\sum_{k=1}^{N-1}
\begin{bmatrix}\Delta x_k^j\\\Delta u_k^j\end{bmatrix}^T\tilde Q_k^j
\begin{bmatrix}\Delta x_k^j\\\Delta u_k^j\end{bmatrix}+
\left(\tilde r_k^j\right)^T\begin{bmatrix}\Delta x_k^j\\\Delta u_k^j\end{bmatrix}\\
&+\frac{1}{2}\left(x_N^j\right)^T\tilde Q_N^jx_N+\left(\tilde r_N^j\right)^T\Delta x_N^j\\
{\rm s.t.}\; &\Delta x_{k+1}^j=A_k^j\Delta x_k^j+B_k^j\Delta u_k^j+\delta r_k^j,\; \Delta x_0=\delta x_0^j\\
&\bar C\Delta u=\delta u
\end{align*}
where $\Delta x_0=(\Delta x_0^1,\ldots,\Delta x_0^M)$,
$\Delta u_0=(\Delta u_0^1,\ldots,\Delta u_0^M)$, and where the other
quantities are defined to agree with the optimality conditions in 
(\ref{eqn:KKT-reduced}). We see that the only 
coupling between the different scenarios are in the first term in the 
objective function and via the non-anticipativity constraints. 
We may equivalently rewrite the above QP  as
\begin{multline}
\min_{\Delta u,\Delta x}\bar F(\Delta x_0,\Delta u_0,)\\
+\sum_{j=1}^M\sum_{k=0}^{N}
\bar F_k^j(\Delta x_k^j,\Delta u_k^j,\Delta x_{k+1}^j,\Delta u_k^{j+1})
\label{prob:QP}
\end{multline}
Here the first function $\bar F$ is the incremental cost for $k=0$.
The remaining functions are the sum of quadratic functions for the 
incremental costs (not for $k=0$)
and indicator 
functions for the constraints, i.e. the dynamic constraints and 
the non-anticipativity constraints. 
We remark that for $k>N_r$
there is no dependence on $\Delta u_k^{j+1}$. Also for smaller values of $k$
this dependence is not present in all $\bar F_k^j$. One has to study the
non-anticipativity constraint in detail to see where it is present. 

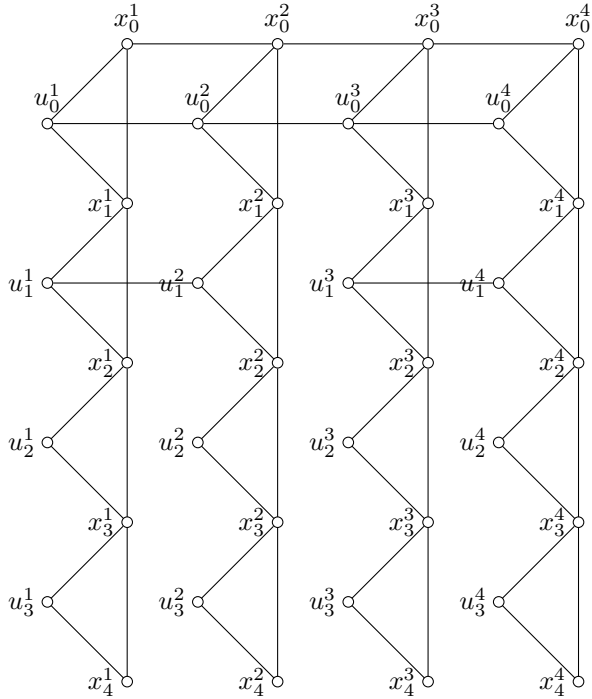
\begin{figure}
		
		\centering
		\begin{tikzpicture}
		\tikzstyle{every node}=[draw,circle,fill=white,minimum size=4pt,
		inner sep=0pt]
		
		\draw (0,0) node (x_01) [label=$x_0^1$] {}
		-- ++(225:1.5cm) node (u_01) [label=$u_0^1$] {}
		-- ++(-45:1.5cm) node (x_11) [label=left:$x_1^1$] {}
		-- (x_01);
		
		\draw (x_11)
		-- ++(225:1.5cm) node (u_11) [label=left:$u_1^1$] {}
		-- ++(-45:1.5cm) node (x_21) [label=left:$x_2^1$] {}
		-- (x_11);
        
        \draw (x_21)
		-- ++(225:1.5cm) node (u_21) [label=left:$u_2^1$] {}
		-- ++(-45:1.5cm) node (x_31) [label=left:$x_3^1$] {}
		-- (x_21);

        \draw (x_31)
		-- ++(225:1.5cm) node (u_31) [label=left:$u_3^1$] {}
		-- ++(-45:1.5cm) node (x_41) [label=left:$x_4^1$] {}
		-- (x_31);
        
        \draw (2,0) node (x_02) [label=$x_0^2$] {}
		-- ++(225:1.5cm) node (u_02) [label=$u_0^2$] {}
		-- ++(-45:1.5cm) node (x_12) [label=left:$x_1^2$] {}
		-- (x_02);
		
		\draw (x_12)
		-- ++(225:1.5cm) node (u_12) [label=left:$u_1^2$] {}
		-- ++(-45:1.5cm) node (x_22) [label=left:$x_2^2$] {}
		-- (x_12);
        
        \draw (x_22)
		-- ++(225:1.5cm) node (u_22) [label=left:$u_2^2$] {}
		-- ++(-45:1.5cm) node (x_32) [label=left:$x_3^2$] {}
		-- (x_22);        

        \draw (x_32)
		-- ++(225:1.5cm) node (u_32) [label=left:$u_3^2$] {}
		-- ++(-45:1.5cm) node (x_42) [label=left:$x_4^2$] {}
		-- (x_32);        
        
        \draw (4,0) node (x_03) [label=$x_0^3$] {}
		-- ++(225:1.5cm) node (u_03) [label=$u_0^3$] {}
		-- ++(-45:1.5cm) node (x_13) [label=left:$x_1^3$] {}
		-- (x_03);
		
		\draw (x_13)
		-- ++(225:1.5cm) node (u_13) [label=left:$u_1^3$] {}
		-- ++(-45:1.5cm) node (x_23) [label=left:$x_2^3$] {}
		-- (x_13);
        
        \draw (x_23)
		-- ++(225:1.5cm) node (u_23) [label=left:$u_2^3$] {}
		-- ++(-45:1.5cm) node (x_33) [label=left:$x_3^3$] {}
		-- (x_23);        
        
        \draw (x_33)
		-- ++(225:1.5cm) node (u_33) [label=left:$u_3^3$] {}
		-- ++(-45:1.5cm) node (x_43) [label=left:$x_4^3$] {}
		-- (x_33);        
        
        \draw (6,0) node (x_04) [label=$x_0^4$] {}
		-- ++(225:1.5cm) node (u_04) [label=$u_0^4$] {}
		-- ++(-45:1.5cm) node (x_14) [label=left:$x_1^4$] {}
		-- (x_04);
		
		\draw (x_14)
		-- ++(225:1.5cm) node (u_14) [label=left:$u_1^4$] {}
		-- ++(-45:1.5cm) node (x_24) [label=left:$x_2^4$] {}
		-- (x_14);
        
        \draw (x_24)
		-- ++(225:1.5cm) node (u_24) [label=left:$u_2^4$] {}
		-- ++(-45:1.5cm) node (x_34) [label=left:$x_3^4$] {}
		-- (x_24);        

        \draw (x_34)
		-- ++(225:1.5cm) node (u_34) [label=left:$u_3^4$] {}
		-- ++(-45:1.5cm) node (x_44) [label=left:$x_4^4$] {}
		-- (x_34);        
        
        \draw (x_01) -- (x_02);
        \draw (x_02) -- (x_03);
        \draw (x_03) -- (x_04);
        
        \draw (u_01) -- (u_02);
        \draw (u_02) -- (u_03);
        \draw (u_03) -- (u_04);

        \draw (u_11) -- (u_12);
        \draw (u_13) -- (u_14);
         
		\end{tikzpicture}
		\caption{Sparsity graph for the problem in 
(\ref{prob:QP}). }
		\label{fig:rob-sparsity}
	\end{figure}

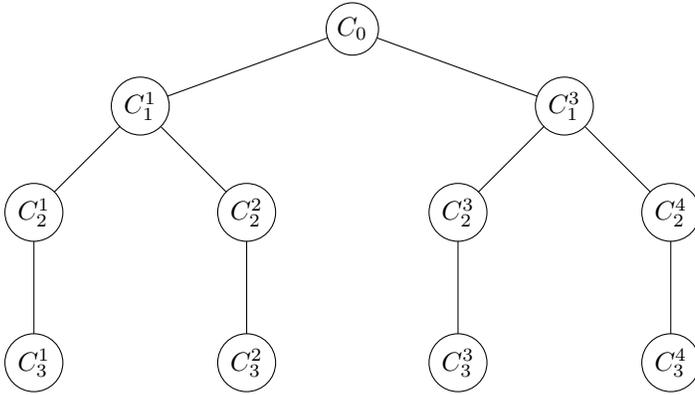
\begin{figure}
		
		\centering
		\begin{tikzpicture}
		
		\tikzstyle{every node}=[draw,circle,fill=white,minimum size=2pt,
		inner sep=2pt]
		
		\draw (0,0) node (0) {$C_0$}
		-- ++(200:3cm) node (11) {$C_1^1$}
		-- ++(225:2cm) node (21) {$C_2^1$}  
		-- ++(270:2cm) node (31) {$C_3^1$}
		;
		
		\draw (0) 
		-- ++(-20:3cm) node (13) {$C_1^3$}  
		-- ++(225:2cm) node (23) {$C_2^3$} 
		-- ++(270:2cm) node (33) {$C_3^3$} 
		;	
        
        \draw (11) 
		-- ++(-45:2cm) node (22) {$C_2^2$}  
		-- ++(270:2cm) node (32) {$C_3^2$} 
		;	
        
        \draw (13) 
		-- ++(-45:2cm) node (24) {$C_2^4$}  
		-- ++(270:2cm) node (34) {$C_3^4$} 
		;	
		\end{tikzpicture}
		
		\caption{Clique trees for the problem in (\ref{prob:QP}).}
		\label{fig:rob-clique}
		
\end{figure}

We now study the case when $N_r=1$, $M_0=M_1=2$ and $N=4$ in more detail. 
Then $M=4$. The sparsity graph is shown in 
Figure~\ref{fig:rob-sparsity}. We label the nodes with $x_0^1$ instead
of $\Delta x_0^1$ and so on. Moreover we do not show all the edges
related to the coupling in $\bar F(\Delta x_0,\Delta u_0)$ since this
would clutter the graph. Actually all of the eight 
variables $x_0^j$ and $u_0^j$ have edges connecting them.
We realize that the sparsity graph is
not chordal. A chordal embedding is obtained by adding
edges such that $C_0=\lbrace x_0^1,_0^2,x_0^3,x_0^4,u_0^1,u_0^2,u_0^3,u_0^4,
x_0^1,_0^2,x_0^3,x_0^4\rbrace$, 
$C_1^1=\lbrace x_0^1,u_0^1,x_1^1,x_0^2,u_0^2,x_1^2\rbrace$ and 
$C_1^3=\lbrace x_0^3,u_0^3,x_1^3,x_0^4,u_0^4,x_1^4\rbrace$ are complete graphs. 
A clique tree for the chordal embedding is 
shown in Figure~\ref{fig:rob-clique},
where $C_{k+1}^j=\lbrace x_k^j,u_k^j,x_{k+1}^j\rbrace$ with $k\in\naturals_{N-1}$.
The assigned functions to $C_0$ are
$$\bar F(\Delta x_0,\Delta u_0,)+
\sum_{j=1}^4\bar F_0^j(\Delta x_0^j,\Delta u_0^j,\Delta x_{1}^j,\Delta u_0^{j+1})
$$
for $C_1^1$
$$
\bar F_0^1(\Delta x_0^1,\Delta u_0^1,\Delta x_{1}^1,\Delta u_1^{2})+
\bar F_0^2(\Delta x_0^2,\Delta u_0^2,\Delta x_{1}^2)
$$
and for $C_1^3$ are
$$
\bar F_0^3(\Delta x_0^3,\Delta u_0^3,\Delta x_{1}^3,\Delta u_1^{4})+
\bar F_0^4(\Delta x_0^4,\Delta u_0^4,\Delta x_{1}4)
$$
For $C_{k+1}^j$, where $k\in\naturals_{N-1}$ and $j\in\naturals_M$, we assign
$$\bar F_{k}^j(\Delta x_k^j,\Delta u_k^j,\Delta x_{k+1}^j)$$
It is possible to introduce even further parallelism by combining the above 
formulation with a parallel formulation in time as described in 
\cite{han+pak17}. It is possible to make use of Riccati recursions to 
compute the messages that are sent up-wards in the clique trees, see
\cite{han+pak17} for details. However, there is no reason to do this. 
A general purpose solver for loosely coupled convex problems with chordal
structure is as efficient and much easier to use. This is the main message
of this article. 
\section{Conclusions}
We have in this paper shown how it is possible to make use of the inherent 
chordal structure of a robust MPC problem in order to exploit IP methods 
that make use of any chordal structure to distribute its computations over
several computational agents that can work in parallel. We argue that this
level of abstraction, i.e. chordality, is more appropriate than a more 
detailed level of abstraction where one tries to see Riccati recursion 
structure. The reason for this is that chordality is a more general concept. 
It also appears when the dynamic equations are obtained from
 spatial discretization of partial differential equations. Hence we believe
that this structure can be utilized using the same formalism as we 
have presented above. How to carry out these extensions is left for future work.
Also it is left for future work to implement a code that carries out the 
computations is parallel and to make comparisons with serial implementations. 
\bibliographystyle{plain}
\bibliography{myref,PaperIII}
\end{document}